\documentclass{article}
\usepackage[dvips]{graphicx}
\usepackage[latin1]{inputenc}
\usepackage{amssymb,amsmath,array}
\usepackage{graphics}
\usepackage[english]{babel}

\begin{document}
\title{Estimating the Number of Components in a Mixture of Multilayer Perceptrons}

\author{M. Olteanu and J. Rynkiewicz
\vspace{.3cm}\\
%
SAMOS-MATISSE-CES Universite Paris 1, UMR 8174 \\
90 Rue de Tolbiac, 75013 Paris, France
}

\maketitle
\begin{abstract}
BIC criterion is widely used by the neural-network community for model selection tasks, although its convergence properties are not always theoretically established. In this paper we will focus on estimating the number of components in a mixture of multilayer perceptrons and proving the convergence of the BIC criterion in this frame. The penalized marginal-likelihood for mixture models and hidden Markov models introduced by Keribin (2000) and, respectively, Gassiat (2002) is extended to mixtures of multilayer perceptrons for which a penalized-likelihood criterion is proposed. We prove its convergence under some hypothesis which involve essentially the bracketing entropy of the generalized score-functions class and illustrate it by some numerical examples.
\end{abstract}

\section{Introduction}

Although linear models have been the standard tool for time series analysis for a long time, their limitations have been underlined during the past twenty years. Real data often exhibit characteristics that are not taken into account by linear models. Financial series, for instance, alternate strong and weak volatility periods, while economic series are often related to the business cycle and switch from recession to normal periods. Several solutions such as heteroscedatic ARCH, GARCH models, threshold models, multilayer perceptrons or autoregressive switching Markov models were proposed to overcome these problems.

In this paper, we consider models which allow the series to switch between regimes and more particularly we study the case of mixtures of multilayer perceptrons. In this frame, rather than using a single global model, we estimate several local models from the data. For the moment, we assume that switches between different models occur independently, the next step of this approach being to also learn how to split the input space and to consider the more general case of \emph{gated experts} or \emph{mixtures of experts} models (Jacobs et al., 1991). The problem we address here is how to select the number of components in a mixture of multilayer perceptrons. This is typically a problem of non-identifiability which leads to a degenerate Fisher information matrix and the classical chi-square theory on the convergence of the likelihood ratio fails to apply. One possible method to answer this problem is to consider penalized criteria. The consistency of the BIC criterion was recently proven for non-identifiable models such as mixtures of densities or hidden Markov models (Keribin, 2000 and Gassiat, 2002). We extend these results to mixtures of nonlinear autoregressive models and prove the consistency of a penalized estimate for the number of components under some good regularity conditions. 

The rest of the paper is organized as follows : in Section 2 we give the definition of the general model and state sufficient conditions for regularity. Then, we introduce the penalized likelihood estimate for the number of components and state the result of consistency. Section 3 is concerned with applying the main result to mixtures of
multilayer perceptrons. Some open questions, as well as some possible extensions are discussed in the conclusion.

\section{Penalized-likelihood estimate for the number of components in a mixture
of nonlinear autoregressive models}

This section is devoted to the setting of the general theoretical frame : model, 
definition and consistency of the penalized-likelihood estimate for the number of components.

\smallskip{}
\textbf{The model - definition and regularity conditions}
\smallskip{}

Throughout the paper, we shall consider that the number of lags is known and, for ease of writing, we shall set the number of lags equal to one, the extension to $l$ time-lags being immediate.

Let us consider the real-valued time series $Y_{t}$ which verifies the following model 

\begin{flushleft}\[
(1)\qquad Y_{t}=F_{\theta_{X_{t}}}^{0}\left(Y_{t-1}\right)+\varepsilon_{\theta_{X_{t}}}\left(t\right)\,,\]
\par\end{flushleft}

where

\begin{itemize}
\item $X_{t}$ is an iid sequence of random variables valued in a finite
space $\left\{ 1,...,p_{0}\right\} $ and with probability distribution
$\pi^{0}$ ;
\item for every $i\in\left\{ 1,...,p_{0}\right\} $, $F_{\theta_{i}}^{0}\left(y\right)\in\mathcal{F}$
and

$\mathcal{F}=\left\{F_{\theta},\,\theta\in\Theta,\,\Theta\subset\mathbb{R}^{l}\,\mathrm{compact}\,\mathrm{set}\right\} $

is the family of possible regression functions. We suppose throughout
the rest of the paper that $F_{\theta_{i}}^{0}$ are sublinear, that
is they are continuous and there exist $\left(a_{i}^{0},b_{i}^{0}\right)\in\mathbb{R}_{+}^{2}$
such that $\left|F_{\theta_{i}}^{0}\left(y\right)\right|\leq a_{i}^{0}\left|y\right|+b_{i}^{0},\:\left(\forall\right)\, y\in\mathbb{R}$
;
\item for every $i\in\left\{ 1,...,p_{0}\right\} $, $\left(\varepsilon_{\theta_{i}}\left(t\right)\right)_{t}$
is an iid noise such that $\varepsilon_{\theta_{i}}\left(t\right)$
is independent of $\left(Y_{t-k}\right)_{k\geq1}$. Moreover, $\varepsilon_{\theta_{i}}\left(t\right)$
has a centered Gaussian density $f_{\theta_{i}}^{0}$. 
\end{itemize}

The sublinearity condition on the regression functions is quite general and the consistency for the number of components  holds for various classes of processes, such as mixtures of densities, mixtures of linear autoregressive functions or mixtures of multilayer perceptrons.

Let us also remark that besides its necessity in the proof of the theoretical result, the compactness hypothesis for the parameter space is also useful in practice. Indeed, one needs to bound the parameter space in order to avoid numerical problems in the multilayer perceptrons such as hidden-unit saturation. In our case, $10^{6}$ seems to be an acceptable bound for the computations. On the other hand, mixture probabilities are naturally bounded.

The next example of a linear mixture illustrates the model introduced by (1). The hidden process $X_{t}$ is a sequence of iid variables with  Bernoulli(0.5) distribution. We define $Y_{t}$ as follows, using $X_{t}$ and a standard Gaussian noise $\varepsilon_{t}$ :

\[Y_{t}=\left\{ \begin{array}{ll} 0.5Y_{t-1}+\varepsilon_{t} & ,if\: X_{t}=1\\
-0.5Y_{t-1}+\varepsilon_{t} & ,if\: X_{t}=0\end{array}\right.\]

The penalized-likelihood estimate which we introduce in the next subsection converges in probability to the true number of components of the model under some regularity conditions on the process $Y_{t}$. More precisely, we need the following hypothesis which implies, according to Yao and Attali (2000), strict stationarity and geometric ergodicity for $Y_{t}$ :

\begin{center}\textbf{(HS)$\quad $} $\sum _{i=1}^{p_{0}}\pi _{i}^{0}\left|a_{i}^{0}\right|^{s}<1$\end{center}

Let us remark that hypothesis \textbf{(HS)} does not request every component to be stationary and that it allows non-stationary {}``regimes'' as long as they do not appear too often. Since multilayer perceptrons are bounded function, this hypothesis will be naturally fulfilled.

\smallskip{}
\textbf{Construction of the penalized likelihood criterion}
\smallskip{}

Let us consider an observed sample $\left\{ y_{1},...,y_{n}\right\} $
of the time series $Y_{t}$. Then, for every observation $y_{t}$,
the conditional density with respect to the previous $y_{t-1}$ and
marginally in $X_{t}$ is

\[
g^{0}\left(y_{t}\mid y_{t-1}\right)=\sum_{i=1}^{p_{0}}\pi_{i}^{0}f_{\theta_{i}}^{0}\left(y_{t}-F_{\theta_{i}}^{0}\left(y_{t-1}\right)\right)\]

As the goal is to estimate $p_{0}$, the number of regimes of the model, let us consider all possible conditional densities up to a maximal number of regimes $P$, a fixed positive integer. We shall consider the class of functions

\[
\mathcal{G}_{P}=\bigcup_{p=1}^{P}\mathcal{G}_{p},\:\mathcal{G}_{p}=\left\{ g\mid g\left(y_{1},y_{2}\right)=\sum_{i=1}^{p}\pi_{i}f_{\theta_{i}}\left(y_{2}-F_{\theta_{i}}\left(y_{1}\right)\right)\right\} ,\]

\medskip{}

where $\pi_{i}\geq\eta>0$, $\sum_{i=1}^{p}\pi_{i}=1$, $F_{\theta_{i}}\left(y\right)\in\mathcal{F}$
and $f_{\theta_{i}}$ is a centered Gaussian density.

For every $g\in\mathcal{G}_{P}$ we define the number of regimes as

\[p\left(g\right)=min\left\{ p\in\left\{ 1,...,P\right\} ,\: g\in\mathcal{G}_{p}\right\} \]

and let $p_{0}=p\left(g^{0}\right)$ be the true number of regimes.

We can now define the estimate $\hat{p}$ as the argument $p\in\left\{ 1,...,P\right\} $
maximizing the penalized criterion

\[(2)\qquad T_{n}\left(p\right)=sup_{g\in\mathcal{G}_{p}}\, l_{n}\left(g\right)-a_{n}\left(p\right)\]

where

\[l_{n}\left(g\right)=\sum_{t=2}^{n}ln\: g\left(y_{t-1},y_{t}\right)\]

is the log-likelihood marginal in $X_{k}$ and $a_{n}\left(p\right)$
is a penalty term.\medskip{}

\smallskip{}
\textbf{Convergence of the penalized likelihood estimate}
\smallskip{}

Several statistical and probabilistic notions such as mixing processes, bracketing entropy or Donsker classes will be used hereafter. For parcimony purposes we shall not remind them, but the reader may refer to Doukhan (1995) and Van der Vaart (2000) for complete monographs on the subject. \\

The consistency of $\widehat{p}$ is given by the next result, which in an extension of Gassiat (2002):

\medskip{}

\textbf{\emph{Theorem 1}} \emph{: Consider the model $\left(Y_{k},\, X_{k}\right)$
defined by (1) and the penalized-likelihood criterion introduced in
(2). Let us introduce the next assumptions :}

\textbf{\emph{$\quad $(A1)}} \emph{$a_{n}\left(\cdot \right)$ is
an increasing function of $p$, $a_{n}\left(p_{1}\right)-a_{n}\left(p_{2}\right)\rightarrow \infty $
when $n\rightarrow \infty $ for every $p_{1}>p_{2}$ and $\frac{a_{n}\left(p\right)}{n}\rightarrow 0$
when $n\rightarrow \infty $ for every $p$}

\textbf{\emph{$\quad $(A2)}} \emph{the model} \emph{$\left(Y_{k},\, X_{k}\right)$
verifies the weak identifiability assumption} \textbf{\emph{(HI)}}

\[\sum _{i=1}^{p}\pi _{i}f_{i}\left(y_{2}-F_{i}\left(y_{1}\right)\right)=\sum _{i=1}^{p_{0}}\pi _{i}^{0}f_{i}^{0}\left(y_{2}-F_{i}^{0}\left(y_{1}\right)\right)\Leftrightarrow \sum _{i=1}^{p}\pi _{i}\delta _{\theta _{i}}=\sum _{i=1}^{p_{0}}\pi _{i}^{0}\delta _{\theta _{i}^{0}}\]

\textbf{\emph{$\quad $(A3)}} \emph{the parameterization $\theta _{i}\rightarrow f_{i}\left(y_{2}-F_{i}\left(y_{1}\right)\right)$
is continuous for every $\left(y_{1},y_{2}\right)$ and there exists
$m\left(y_{1},y_{2}\right)$ an integrable map with respect to the
stationary measure of $\left(Y_{k},Y_{k-1}\right)$ such that $\left|log\left(g\right)\right|<m$}

\textbf{\emph{$\quad $(A4)}} \emph{$Y_{k}$ is strictly stationary
and geometrically $\beta $-mixing, and the family of generalized
score functions associated to $\mathcal{G}_{P}$ }

\begin{center}
$\mathcal{S}=\left\{ s_{g},\: s_{g}\left(y_{1},y_{2}\right)=\frac{\frac{g\left(y_{1},y_{2}\right)}{f\left(y_{1},y_{2}\right)}-1}{\left\Vert \frac{g}{f}-1\right\Vert _{L^{2}\left(\mu \right)}},\, g\in \mathcal{G}_{P},\, \left\Vert \frac{g}{f}-1\right\Vert _{L^{2}\left(\mu \right)}\neq 0\right\} \subset \mathcal{L}_{2}\left(\mu \right)$
\end{center}

\emph{where $\mu $ is the stationary measure of} $\left(Y_{k},Y_{k-1}\right)$
\emph{and for every $\varepsilon >0$}

\begin{center}\emph{$\mathcal{H}_{\left[\cdot \right]}\left(\varepsilon ,\mathcal{S},\left\Vert \cdot \right\Vert _{2}\right)=\mathcal{O}\left(\left|log\, \varepsilon \right|\right),$}\end{center}

\emph{$\mathcal{H}_{\left[\cdot \right]}\left(\varepsilon ,\mathcal{S},\left\Vert \cdot \right\Vert _{2}\right)$
being the bracketing entropy of $\mathcal{S}$ with respect to the $L_{2}$-norm. }

\emph{Then, under hypothesis (A1)-(A4), (HS) et (HC), $\hat{p}\rightarrow p_{0}$ in probability.}

\bigskip{}

\textit\textbf{Proof of Theorem 1}

First, let us show that $\widehat{p}$ does not overestimate $p_{0}$. We shall need the following likelihood ratio inequality  which is an immediate generalization of Gassiat (2002) to multivariate dependent data.

\medskip{}

\emph{
Let $\mathcal{G}\subset\mathcal{G}_{P}$ be a parametric family
of conditional densities containing the true model $g^{0}$ and consider
the generalized score functions }

\[s_{g}\left(y_{1},y_{2}\right)=\frac{\frac{g\left(y_{1},y_{2}\right)}{g^{0}\left(y_{1},y_{2}\right)}-1}{\left\Vert \frac{g}{g^{0}}-1\right\Vert _{L^{2}\left(\mu\right)}}\]

\emph{
where $\mu$ is the stationary measure of $\left(Y_{k-1},Y_{k}\right)$.
Then, }

\[
sup_{g\in\mathcal{G}}\,\left(l_{n}\left(g\right)-l_{n}\left(g^{0}\right)\right)\leq\frac{1}{2}sup_{g\in\mathcal{G}}\,\frac{\left(\sum_{k=2}^{n}s_{g}\left(y_{k-1},y_{k}\right)\right)^{2}}{\sum_{k=2}^{n}\left(s_{g}\right)_{\_}^{2}\left(y_{k-1},y_{k}\right)},\]

\emph{with}

 $\left(s_{g}\right)_{\_}\left(y_{k-1},y_{k}\right)=min\left(0,\, s_{g}\left(y_{k-1},y_{k}\right)\right)$.

\bigskip{}

Then we have :

\[
\mathbb{P}\left(\hat{p}>p_{0}\right)\leq\sum_{p=p_{0}+1}^{P}\mathbb{P}\left(T_{n}\left(p\right)>T_{n}\left(p_{0}\right)\right)=\]

\[
=\sum_{p=p_{0}+1}^{P}\mathbb{P}\left(sup_{g\in\mathcal{G}_{p}}l_{n}\left(g\right)-a_{n}\left(p\right)>sup_{g\in\mathcal{G}_{p_{0}}}l_{n}\left(g\right)-a_{n}\left(p_{0}\right)\right)\leq\]

\[
\leq\sum_{p=p_{0}+1}^{P}\mathbb{P}\left(\frac{1}{2}sup_{g\in\mathcal{G}_{p}}\frac{\left(\sum_{k=2}^{n}s_{g}\left(Y_{k-1},Y_{k}\right)\right)^{2}}{\sum_{k=2}^{n}\left(s_{g}\right)_{\_}^{2}\left(Y_{k-1},Y_{k}\right)}>a_{n}\left(p\right)-a_{n}\left(p_{0}\right)\right)\]

Under the hypothesis \textbf{(HS)}, there exists a unique strictly stationary solution $Y_{k}$ which is also geometrically ergodic and this implies that $Y_{k}$ is in particular geometrically $\beta$-mixing. Then, by remarking that

\[
\beta_{n}^{\left(Y_{k-1},Y_{k}\right)}=\beta_{n-1}^{Y_{k}}\]
we obtain that the bivariate series $\left(Y_{k-1},Y_{k}\right)$
is also strictly stationary and geometrically $\beta$-mixing. 

This fact, together with the assumption on the $\varepsilon$-bracketing entropy of $\mathcal{S}$ with respect to the $\left\Vert \cdot\right\Vert _{L^{2}\left(\mu\right)}$ norm and the condition that $\mathcal{S}\subset\mathcal{L}_{2}\left(\mu\right)$ ensures that Theorem 4 in Doukan, Massart and Rio (1995) holds and 

\[\left\{ \frac{1}{\sqrt{n-1}}\sum_{k=2}^{n}s_{g}\left(Y_{k-1},Y_{k}\right)\,\mid g\in\mathcal{G}_{p}\right\} \]

\bigskip{}

is uniformly tight and verifies a functional central limit theorem. Then,

\[
sup_{g\in\mathcal{G}_{p}}\frac{1}{n-1}\left(\sum_{k=2}^{n}s_{g}\left(Y_{k-1},Y_{k}\right)\right)^{2}=\mathcal{O}_{\mathbb{P}}\left(1\right)\]

On the other hand, $\mathcal{S}\subset\mathcal{L}_{2}\left(\mu\right)$,
thus $\mathcal{S}^{2}\subset\mathcal{L}_{1}\left(\mu\right)$ and
using the $\mathcal{L}_{2}$-entropy condition $\mathcal{S}_{\_}^{2}=\left\{ \left(s_{g}\right)_{\_}^{2},\, g\in\mathcal{G}_{p}\right\} $
is Glivenko-Cantelli. Since $\left(Y_{k-1},Y_{k}\right)$ is ergodic
and strictly stationary, we obtain the following uniform convergence
in probability :

\[
inf_{g\in\mathcal{G}_{p}}\frac{1}{n-1}\sum_{k=2}^{n}\left(s_{g}\right)_{\_}^{2}\left(Y_{k-1},Y_{k}\right)\longrightarrow_{n\rightarrow\infty}inf_{g\in\mathcal{G}_{p}}\left\Vert \left(s_{g}\right)_{\_}\right\Vert _{2}^{2}\]
\bigskip{}

To finish the first part, let us prove that

\[inf_{g\in\mathcal{G}_{p}}\left\Vert \left(s_{g}\right)_{\_}\right\Vert _{2}>0\]

If we suppose, on the contrary, that $inf_{g\in\mathcal{G}_{p}}\left\Vert \left(s_{g}\right)_{\_}\right\Vert _{2}=0$,
then there exists a sequence of functions $\left(s_{g_{n}}\right)_{n\geq1}$
, $g_{n}\in\mathcal{G}_{p}$ such that $\left\Vert \left(s_{g_{n}}\right)_{\_}\right\Vert _{2}\rightarrow0$.
The $L_{2}$-convergence implies that $\left(s_{g_{n}}\right)_{\_}\rightarrow0$
in $L_{1}$ and a.s. for a subsequence $s_{g_{n,k}}$. Since $\int s_{g_{n}}d\mu=0$
and $s_{g_{n}}=\left(s_{g_{n}}\right)_{\_}+\left(s_{g_{n}}\right)_{+}$,
where $\left(s_{g_{n}}\right)_{+}=max\left(0,\, s_{g_{n}}\right)$,
we obtain that $\int\left(s_{g_{n}}\right)_{+}d\mu=-\int\left(s_{g_{n}}\right)_{-}d\mu=\int\left|\left(s_{g_{n}}\right)_{-}\right|d\mu$ and thus $\left(s_{g_{n}}\right)_{+}\rightarrow0$ in $L_{1}$ and a.s. for a subsequence $s_{g_{n,k'}}$. The hypothesis \textbf{(A4)} ensures the existence of a square-integrable dominating function for $\mathcal{S}$ and, finally, we get that a subsequence of $s_{g_{n}}$ converges to $0$ a.s. and in $L_{2}$, which contradicts the fact that $\int s_{g}^{2}d\mu=1$ for every $g\in\mathcal{G}_{p}$, so that :

\[
sup_{g\in\mathcal{G}_{p}}\frac{\left(\sum_{k=2}^{n}s_{g}\left(Y_{k-1},Y_{k}\right)\right)^{2}}{\sum_{k=2}^{n}\left(s_{g}\right)_{\_}^{2}\left(Y_{k-1},Y_{k}\right)}=\mathcal{O}_{\mathbb{P}}\left(1\right)\]

Then, by the uniform tightness above and the hypothesis \textbf{(A1)},

\[\mathbb{P}\left(\hat{p}>p_{0}\right)\longrightarrow_{n\rightarrow\infty}0\]

Let us now prove that $\hat{p}$ does not underestimate $p_{0}$ :

\[
\mathbb{P}\left(\hat{p}<p_{0}\right)\leq\sum_{p=1}^{p_{0}-1}\mathbb{P}\left(T_{n}\left(p\right)>T_{n}\left(p_{0}\right)\right)\leq\]

\[
\leq\sum_{p=1}^{p_{0}-1}\mathbb{P}\left(\frac{sup_{g\in\mathcal{G}_{p}}\left(l_{n}\left(g\right)-l_{n}\left(g^{0}\right)\right)}{n-1}>\frac{a_{n}\left(p\right)-a_{n}\left(p_{0}\right)}{n-1}\right)\]

Now, $l_{n}\left(g\right)-l_{n}\left(g^{0}\right)=\sum_{k=2}^{n}ln\left(\frac{g\left(Y_{k-1},Y_{k}\right)}{g^{0}\left(Y_{k-1},Y_{k}\right)}\right)$
and under the hypothesis \textbf{(A3)}, the class of functions $\left\{ ln\frac{g}{g^{0}},\, g\in\mathcal{G}_{p}\right\} $
is $\mathbb{P}$-Glivenko-Cantelli (the general proof for a parametric
family can be found in Van der Vaart, 2000) and since $\left(Y_{k-1},Y_{k}\right)$
is ergodic and strictly stationary, we obtain the following uniform
convergence in probability :

\[
\frac{1}{n-1}sup_{g\in\mathcal{G}_{p}}\left(l_{n}\left(g\right)-l_{n}\left(g^{0}\right)\right)\longrightarrow sup_{g\in\mathcal{G}_{p}}\int ln\frac{g}{g^{0}}g^{0}d\mu\]

Since $p<p_{0}$ and using assumption \textbf{(A2)}, the limit is
negative. By hypothesis \textbf{(A1)}, $\frac{a_{n}\left(p\right)-a_{n}\left(p_{0}\right)}{n-1}$
converges to $0$ when $n\rightarrow\infty$, so we finally have that
$\mathbb{P}\left(\hat{p}<p_{0}\right)\rightarrow0$ and the proof
is done.\bigskip{}

\begin{flushright}$\blacksquare$\par\end{flushright}

\section{Mixtures of multilayer perceptrons}

In this section, we consider the model defined in (1) such that, for every $i\in \left\{ 1,...,p_{0}\right\} $, $F_{i}^{0}$ is a multilayer perceptron. Since non-identifiability problems also arise in multilayer perceptrons (see, for instance, Rynkiewicz, 2006), we shall simplify the problem by considering one hidden layer and a fixed number of
units on every layer, $k$. Then, we have that for every $i\in \left\{ 1,...,p_{0}\right\} $

\begin{center}$F_{i}^{0}\left(y\right)=\alpha _{0}^{0,i}+\sum _{j=1}^{k}\alpha _{j}^{0,i}\phi \left(\beta _{0,j}^{0,i}+\beta _{1,j}^{0,i}y\right)$\end{center}

where $\phi $ is the hyperbolic tangent and

\begin{center}$\theta _{i}^{0}=\left(\alpha _{0}^{0,i},\alpha _{1}^{0,i},...,\alpha _{k}^{0,i},\beta _{0,1}^{0,i},\beta _{1,1}^{0,i},...,\beta _{0,k}^{0,i},\beta _{1,k}^{0,i},\sigma^{0,i} \right)$ \end{center}

is the true parameter with the true variance.Let us check if the hypothesis of the main result of section 2 apply
in the case of mixtures of multilayer perceptrons.

\textbf{Hypothesis (HS)} : The stationarity and ergodicity assumption (HS) is immediately verified since the output of every perceptron is bounded, by construction. Thus, every regime is stationary and the global model is also stationary.

Let us consider the class of all possible conditional densities up to a maximum number of components $P>0$ : 

\begin{center}
$\mathcal{G}_{P}=\bigcup _{p=1}^{P}\mathcal{G}_{p}\: ,\: \mathcal{G}_{p}=\left\{ g\mid g\left(y_{1},y_{2}\right)=\sum _{i=1}^{p}\pi _{i}f_{i}\left(y_{2}-F_{i}\left(y_{1}\right)\right)\right\}$, where
\end{center}

\begin{itemize}
\item $\sum _{i=1}^{p}\pi _{i}=1$ and we may suppose quite naturally that
for every $i\in \left\{ 1,...,p\right\} $, $\pi _{i}\geq \eta >0$ 
\item for every $i\in \left\{ 1,...,p\right\} $, $F_{i}$ is a multilayer
perceptron
\end{itemize}
\begin{center}$F_{i}\left(y\right)=\alpha _{0}^{i}+\sum _{j=1}^{k}\alpha _{j}^{i}\phi \left(\beta _{0,j}^{i}+\beta _{1,j}^{i}y\right)$, where \end{center}

$\theta _{i}=\left(\alpha _{0}^{i},\alpha _{1}^{i},...,\alpha _{k}^{i},\beta _{0,1}^{i},\beta _{1,1}^{i},...,\beta _{0,k}^{i},\beta _{1,k}^{i},\sigma^{i} \right)$
belongs to a compact set.

\textbf{Hypothesis (A1)} : $a_{n}\left(\cdot \right)$ may be chosen, for instance, equal to
the BIC penalizing term, $a_{n}\left(p\right)=\frac{1}{2}p\, log\left(n\right)$. 

\textbf{Hypothesis (A2)-(A3)} : Since the noise is normally distributed, the weak identifiability hypothesis is verified according to the result of Teicher (1963), while assumption (A3) is a regularity condition verified by Gaussian densities.

\textbf{Hypothesis (A4)} : We consider the class of generalized score functions

\begin{center}
$\mathcal{S}=\left\{ s_{g},\: s_{g}=\frac{\frac{g}{f}-1}{\left\Vert \frac{g}{f}-1\right\Vert _{L^{2}\left(\mu \right)}},\: g\in \mathcal{G}_{P},\, \left\Vert \frac{g}{f}-1\right\Vert _{L^{2}\left(\mu \right)}\neq 0\right\} $
\end{center}

The difficult part will be to show that $\mathcal{H}_{\left[\cdot \right]}\left(\varepsilon ,\mathcal{S},\left\Vert \cdot \right\Vert _{2}\right)=\mathcal{O}\left(\left|log\, \varepsilon \right|\right)$ for all $\varepsilon >0$ which, since we are on a functional space,
is equivalent to prove that {}``the dimension'' of $\mathcal{S}$
can be controlled. For $g\in \mathcal{G}_{p}$, let us denote $\theta =\left(\theta _{1},...,\theta _{p}\right)$
and $\pi =\left(\pi _{1},...,\pi _{p}\right)$, so that the global
parameter will be $\Phi =\left(\theta ,\pi \right)$ and the associated
generalized score function $s_{\Phi }:=s_{g}$.

Proving that a parametric family like $\mathcal{S}$ verifies the condition on the bracketing entropy is usually immediate under good regularity conditions (see, for instance, Van der Vaart, 2000). A sufficient condition is that the bracketing number grows as a polynomial function of $\frac{1}{\varepsilon}$.  In this particular case, the problems arise when $g\rightarrow f$ and the limits in $L^{2}\left(\mu \right)$ of $s_{g}$ have to be computed. Let us split $\mathcal{S}$ into two classes of functions. We shall consider $\mathcal{F}_{0}\subset \mathcal{G}_{P}$ a neighborhood of $f$ such that 
$\mathcal{F}_{0}=\left\{ g\in G_{p},\: \left\Vert \frac{g}{f}-1\right\Vert _{L^{2}\left(\mu \right)}\leq \varepsilon\right\} $
and $\mathcal{S}_{0}=\left\{ s_{g},\: g\in \mathcal{F}_{0}\right\} $. On $\mathcal{S}\setminus \mathcal{S}_{0}$, it can be easily seen that 

\begin{center}$\left\Vert \frac{\frac{g_{1}}{f}-1}{\left\Vert \frac{g_{1}}{f}-1\right\Vert _{L^{2}\left(\mu \right)}}-\frac{\frac{g_{2}}{f}-1}{\left\Vert \frac{g_{2}}{f}-1\right\Vert _{L^{2}\left(\mu \right)}}\right\Vert _{L^{2}\left(\mu \right)}\leq \frac{2}{\varepsilon }\left\Vert \frac{g_{1}}{f}-\frac{g_{2}}{f}\right\Vert _{L^{2}\left(\mu \right)}$\end{center}
Hence, on $\mathcal{S}\setminus \mathcal{S}_{0}$, it is sufficient that
\[\left\Vert\frac{g_{1}}{f}-\frac{g_{2}}{f}\right\Vert_2<\frac{\varepsilon^2}{2}\]
for
\[
\left\Vert\frac{\frac{g_{1}}{f}-1}{\left\Vert\frac{g_{1}}{f}-1\right\Vert_2}-\frac{\frac{g_{_2}}{f}-1}{\left\Vert\frac{g_{2}}{f}-1\right\Vert_2}\right\Vert_2<\varepsilon.\]

Now, $\mathcal{S}\setminus \mathcal{S}_{0}$ is a parametric class. Since the derivatives of the transfer functions are bounded,
according to the example 19.7 of Van der Vaart (2000), it exists a constant $K$ so that the bracketing number of ${\mathbb S}_\varepsilon$ is lower than
\[
K\left(\frac{\mbox{diam}G_p}{\varepsilon^2}\right)^{3(k+1)P}=K\left(\frac{\sqrt{\mbox{diam}G_p}}{\varepsilon}\right)^{6(k+1)P},\]
where $\mbox{diam}G_p$ is the diameter of the smallest sphere of $\mathbb R^{3(k+1)P}$ including the set of possible parameters. 
So, we get that $\mathcal{N}_{\left[\right]}\left(\varepsilon ,\mathcal{S}\setminus \mathcal{S}_{0},\left\Vert \cdot \right\Vert _{2}\right)=\mathcal{O}\left(\frac{1}{\varepsilon }\right)^{6(k+1)P}$,
where $\mathcal{N}_{\left[\right]}\left(\varepsilon ,\mathcal{S}\setminus \mathcal{S}_{0},\left\Vert \cdot \right\Vert _{2}\right)$
is the number of $\varepsilon $-brackets necessary to cover $\mathcal{S}\setminus \mathcal{S}_{0}$
and the bracketing entropy is computed as $log\mathcal{N}_{\left[\right]}\left(\varepsilon ,\mathcal{S}\setminus \mathcal{S}_{0},\left\Vert \cdot \right\Vert _{2}\right)$.

As for $\mathcal{S}_{0}$, the idea is to reparameterize the model
in a convenient manner which will allow a Taylor expansion around
the identifiable part of the true value. For that, we shall use a
slight modification of the method proposed by Liu and Shao (2003). Let us remark that
when $\frac{g}{f}-1=0$, the weak identifiability
hypothesis \textbf{(A2)} and the fact that for every $i\in \left\{ 1,...,p\right\} $,
$\pi _{i}\geq \eta >0$, imply that there exists a vector $t=\left(t_{i}\right)_{0\leq i\leq p_{0}}$
such that $0=t_{0}<t_{1}<...<t_{p_{0}}=p$ and, modulo a permutation,
$\Phi $ can be rewritten as follows : $\theta _{t_{i-1}+1}=...=\theta _{t_{i}}=\theta _{i}^{0}$, $\sum _{j=t_{i-1}+1}^{t_{i}}\pi _{j}=\pi _{i}^{0}$, $i\in \left\{ 1,...,p_{0}\right\} $.
With this remark, one can define in the general case $s=\left(s_{i}\right)_{1\leq i\leq p_{0}}$
and $q=\left(q_{j}\right)_{1\leq j\leq p}$ so that, for every $i\in \left\{ 1,...,p_{0}\right\} $
, $j\in \left\{ t_{i-1}+1,...,t_{i}\right\} $,

\begin{center}$s_{i}=\sum _{j=t_{i-1}+1}^{t_{i}}\pi _{j}-\pi _{i}^{0},\; q_{j}=\frac{\pi _{j}}{\sum _{l=t_{i-1}+1}^{t_{i}}\pi _{l}}$\end{center}

and the new parameterization will be $\Theta _{t}=\left(\phi _{t},\psi _{t}\right)$,

$\phi _{t}=\left(\left(\theta _{j}\right)_{1\leq j\leq p},\left(s_{i}\right)_{1\leq i\leq p_{0}-1}\right)$, $\psi _{t}=\left(q_{j}\right)_{1\leq j\leq p}$,
with $\phi _{t}$ containing all the identifiable parameters of the
model and $\psi _{t}$ the non-identifiable ones. Then, for $g=f$,
we will have that

\begin{center}$\phi _{t}^{0}=\begin{array}[t]{cccc}
 (\underbrace{\theta _{1}^{0},...,\theta _{1}^{0}} & ,..., & \underbrace{\theta _{p_{0}}^{0},...,\theta _{p_{0}}^{0}}, & \underbrace{0,...,0}\\
 t_{1} &  & t_{p_{0}}-t_{p_{0}-1} & p_{0}-1\end{array})^{T}$\end{center}

This reparameterization allows to write a second-order Taylor expansion
of $\frac{g}{f}-1$ at $\phi _{t}^{0}$. For ease of writing, we shall
first denote

\begin{center}$g_{j}\left(y_{1},y_{2}\right)=g_{\theta _{j}}\left(y_{1},y_{2}\right)=\frac{f_{j}\left(y_{2}-F_{j}\left(y_{1}\right)\right)}{\sum _{i=1}^{p_{0}}\pi _{i}^{0}f_{i}^{0}\left(y_{2}-F_{i}^{0}\left(y_{1}\right)\right)}-1$\end{center}

Then, the density ratio becomes :

\begin{center}$\frac{g}{f}-1=\sum _{i=1}^{p_{0}-1}\left(s_{i}+\pi _{i}^{0}\right)\sum _{j=t_{i-1}+1}^{t_{i}}q_{j}g_{j}+\left(\pi _{p_{0}}^{0}-\sum _{i=1}^{p_{0}-1}s_{i}\right)\sum _{j=t_{p_{0}-1}+1}^{t_{p_{0}}}q_{j}g_{j}$\end{center}

By remarking that when $\phi _{t}=\phi _{t}^{0}$, $\frac{g}{f}$
does not vary with $\psi _{t}$, we will study the variation of this
ratio in a neighborhood of $\phi _{t}^{0}$ and for fixed $\psi _{t}$.
Let us note $\frac{\partial g_{j}}{\partial \theta _{j}}$ the vector of derivatives of $g_j$ with respect of each components of $\theta _{j}$ and 
$\frac{\partial ^{2}g_{j}}{\partial \theta _{j}^{2}}$ the vector of second derivatives of $g_j$ with respect of each components of $\theta _{j}$.
Assuming that $\left(g_{j}\right)_{1\leq j\leq p}$, $\left(g_{j}^{\prime }\right)_{1\leq j\leq p}$ and
$\left(g_{j}^{\prime \prime }\right)_{1\leq j\leq p}$, where
\begin{center}
$g_{j}^{\prime }:=\frac{\partial g_{j}}{\partial \theta _{j}}\left(\phi _{t}^{0},\psi _{t}\right),\: g_{j}^{\prime \prime }:=\frac{\partial ^{2}g_{j}}{\partial \theta _{j}^{2}}\left(\phi _{t}^{0},\psi _{t}\right)$
\end{center}
are linearly independent in $L^{2}\left(\mu \right)$,
one can prove the following :

\textbf{\emph{Proposition 1}} \emph{: Let us denote $D\left(\phi _{t},\psi _{t}\right)=\left\Vert \frac{g_{\left(\phi _{t},\psi _{t}\right)}}{f}-1\right\Vert _{L^{2}\left(\mu \right)}$.
For any fixed $\psi _{t}$, there exists the second-order Taylor expansion
at $\phi _{t}^{0}$ :}

\begin{center}
$\frac{g}{f}-1=\left(\phi _{t}-\phi _{t}^{0}\right)^{T}g_{\left(\phi _{t}^{0},\psi _{t}\right)}^{\prime }+\frac{1}{2}\left(\phi _{t}-\phi _{t}^{0}\right)^{T}g_{\left(\phi _{t}^{0},\psi _{t}\right)}^{\prime \prime }\left(\phi _{t}-\phi _{t}^{0}\right)+o\left(D\left(\phi _{t},\psi _{t}\right)\right),$
\end{center}

\emph{with}

\[
\left(\phi_{t}-\phi_{t}^{0}\right)^{T}g_{\left(\phi_{t}^{0},\psi_{t}\right)}^{\prime}=\sum_{i=1}^{p_{0}}\pi_{i}^{0}\left(\sum_{j=t_{i-1}+1}^{t_{i}}q_{j}\theta_{j}-\theta_{i}^{0}\right)^{T}g_{i}^{\prime}+\sum_{i=1}^{p_{0}}s_{i}g_{\theta_{i}^{0}}\]

\emph{and}

\[\left(\phi_{t}-\phi_{t}^{0}\right)^{T}g_{\left(\phi_{t}^{0},\psi_{t}\right)}^{\prime\prime}\left(\phi_{t}-\phi_{t}^{0}\right)=\sum_{i=1}^{p_{0}}\left[2s_{i}\left(\sum_{j=t_{i-1}+1}^{t_{i}}q_{j}\theta_{j}-\theta_{i}^{0}\right)^{T}g_{i}^{\prime}+\right.\]

\[\left.+\pi_{i}^{0}\sum_{j=t_{i-1}+1}^{t_{i}}q_{j}\left(\theta_{j}-\theta_{i}^{0}\right)^{T}g_{i}^{\prime\prime}\left(\theta_{j}-\theta_{i}^{0}\right)\right]\]

\emph{Moreover,}
\begin{center}
$\left(\phi _{t}-\phi _{t}^{0}\right)^{T}g_{\left(\phi _{t}^{0},\psi _{t}\right)}^{\prime }+\frac{1}{2}\left(\phi _{t}-\phi _{t}^{0}\right)^{T}g_{\left(\phi _{t}^{0},\psi _{t}\right)}^{\prime \prime }\left(\phi _{t}-\phi _{t}^{0}\right)=0\Leftrightarrow \phi _{t}=\phi _{t}^{0}$
\end{center}

\medskip{}

\textbf{\emph{Proof of Proposition 1}}\medskip{}

The first term in the developpement can be computed easily by remarking
that the gradient of $\frac{g}{g^{0}}-1$ at $\left(\phi_{t}^{0},\psi_{t}\right)$
is :

\begin{itemize}
\item for $i\in\left\{ 1,...,p_{0}\right\} $ and $j\in\left\{ t_{i-1}+1,...,t_{i}\right\} $,
$\frac{\partial\left(\frac{g}{g^{0}}-1\right)}{\partial\theta_{j}}\left(\phi_{t}^{0},\psi_{t}\right)=\pi_{i}^{0}q_{j}g_{i}^{\prime}$
\item for $i\in\left\{ 1,...,p_{0}-1\right\} $,
\end{itemize}
\begin{center}$\frac{\partial\left(\frac{g}{g^{0}}-1\right)}{\partial s_{i}}\left(\phi_{t}^{0},\psi_{t}\right)=\sum_{j=t_{i-1}+1}^{t_{i}}q_{j}g_{\theta_{i}^{0}}-\sum_{j=t_{p_{0}-1}+1}^{t_{p_{0}}}q_{j}g_{\theta_{p_{0}}^{0}}=g_{\theta_{i}^{0}}-g_{\theta_{p_{0}}^{0}}$\par\end{center}

The term of second order can be obtained by direct computations once
the hessian in computed at $\left(\phi_{t}^{0},\psi_{t}\right)$:

\begin{itemize}
\item $\frac{\partial^{2}\left(\frac{g}{g^{0}}-1\right)}{\partial\theta_{j}^{2}}\left(\phi_{t}^{0},\psi_{t}\right)=\pi_{i}^{0}q_{j}g_{i}^{\prime\prime}$
, $i=1,...,p_{0}$ and $j=t_{i-1}+1,...,t_{i}$
\item $\frac{\partial^{2}\left(\frac{g}{g^{0}}-1\right)}{\partial\theta_{j}\partial\theta_{l}}\left(\phi_{t}^{0},\psi_{t}\right)=0$
, $j,l=1,...,p$ and $j\neq l$
\item $\frac{\partial^{2}\left(\frac{g}{g^{0}}-1\right)}{\partial s_{i}\partial s_{k}}\left(\phi_{t}^{0},\psi_{t}\right)=0$
, $i,k=1,...,p_{0}-1$
\item $\frac{\partial^{2}\left(\frac{g}{g^{0}}-1\right)}{\partial s_{i}\partial\theta_{j}}\left(\phi_{t}^{0},\psi_{t}\right)=q_{j}g_{i}^{\prime}$
, $i=1,...,p_{0}-1$ and $j=t_{i-1}+1,...,t_{i}$
\item $\frac{\partial^{2}\left(\frac{g}{g^{0}}-1\right)}{\partial s_{i}\partial\theta_{j}}\left(\phi_{t}^{0},\psi_{t}\right)=-q_{j}g_{p_{0}}^{\prime}$
, $i=1,...,p_{0}-1$ and $j=t_{p_{0}-1}+1,...,t_{p_{0}}$
\item the other crossed derivatives of $s_{i}$ and $\theta_{j}$ are zero
\end{itemize}
It remains to prove that the rest is $o\left(\left\Vert \phi_{t}-\phi_{t}^{0}\right\Vert \right)$ but this follows directly from Yao (2000) and the fact that, since the noise is normally distributed, $Y_{t}$ has moments of any order.

\begin{flushright}$\blacksquare$\par\end{flushright}

\bigskip{}

Using the Taylor expansion above, for $\theta$ belonging to ${\mathbb S}_0$, \(\frac{f_\theta(z)}{f}-1\) is the sum of a linear combination of
\[V(z):=\left(g_1,\cdots,g_{p},g^{'}_1,\cdots,g^{'}_{p},g^{''}_1,\cdots,g^{''}_{p}\right)\] 

and of a term whose $L^2$ norm is negligible compared to the $L^2$ norm of this combination when $\varepsilon$ goes to 0. 
By assumption (A3), a strictly positive number $m$ exists so that for any vector of norm 1 with components
\[C=\left(c_1,\cdots,c_{p_0\times(3k+1)},d_1,\cdots,d_{p_0\times(3k+1)},e_1,\cdots,e_{p_0\times(3k+1)}\right)\]
and $\varepsilon$ sufficiently small:
\[\Vert C^TV(z)\Vert_2> m+\varepsilon.\]
Since any function $\frac{\frac{g}{f}-1}{\Vert \frac{g}{f}-1\Vert_2}$ can be written:
\[\frac{C^TV(z)+o(\Vert C^TV(z)\Vert_2)}{\Vert C^TV(z)+o(\Vert C^TV(z)\Vert_2)\Vert_2},\] 
$\mathbb S_0$ belongs to the set of functions:
\[\left\{D^TV(z)+o(1),\Vert D\Vert_2\leq\frac{1}{m}\right\}\subset\left\{D^TV(z)+\gamma,\Vert D\Vert_2\leq\frac{1}{m},|\gamma|<1\right\}\] 
whose bracketing number is smaller or equal to $O\left(\frac{1}{\varepsilon}\right)^{3p_0\times(3k+1)+1}$.

and the assumptions of Theorem 1 are verified $\blacksquare$

\section{Some numerical examples}

The theoretical result proven above may be applied in practice to compute the number of components in a mixture model on simulated or real data. Some examples are presented below, illustrating the stability  and the speed of convergence of the algorithm. 

\subsection{Mixtures of linear models}

Let us first consider the particular case of linear models, corresponding to zero hidden-unit perceptrons. The examples are mixtures of two autoregressive models in which we vary the leading coefficients and the weights of the discrete mixing distribution. For every example, we simulate 20 samples of lengths $n=200,500,1000,1500,2000$ and we fix $P=3$ the upper bound for the number of regimes.

The likelihood is maximized via the EM algorithm (see, for instance,
Dempster, Laird and Rubin (1977) or Redner and Walker, 1984). This algorithm is well suited to find a sequence of parameters which increases the likelihood at each step, and so converges to a local maximum for a very wide class of models and for our model in particular. The idea of the EM algorithm is to replace the latent variables of the mixture by their conditional expectation. A brief recall on the main steps of the algorithm is given below :

\begin{enumerate}
\item Let $X$ be the vector containing the component of the mixture and considered as a latent variable and let  $y=(y_1,\cdots,y_n)$ be the vector of observations.
\item Initialization : Set \( k=0 \) and choose \( \theta _{0} \) 
\item E-Step : Set \( \theta ^{*}=\theta _{k} \) and compute \( Q(.,\theta ^{*}) \)
with \\
\( Q(\theta ,\theta ^{*})=E_{\theta ^{*}}\left[ \ln \left( \frac{L_{\theta }(y,X)}{L_{\theta ^{*}}(y,X)}\right) \right]  \) where $L_{\theta }(y,X)$ is the likelihood of the observations and the vector of mixture $X$ for the parameter $\theta$. This step computes the probabilities of the mixture conditionally to the observations and with respect to the parameter  $\theta^*$
\item M-Step : Find :\\
\( \hat{\theta }=\arg \max Q(\theta ,\theta ^{*}) \)
\item Replace \( \theta _{k+1} \) by \( \hat{\theta } \) , and go back to
 step 3) until a stopping criterion is satisfied (i.e. when the parameters don't seem to change anymore).
\end{enumerate}
The sequence \( (\theta _{k}) \) gives nondecreasing values of the likelihood
function up to a local maximum. \( Q(\theta ,\theta ^{*}) \) is called conditional pseudo-log-likelihood.

\bigskip{}

To avoid local maxima, the procedure is initialized several times with different starting values : in our case, ten different initializations provided good results. The stopping criteria applies when either there is no improvement in the likelihood value, either a maximum number of iterations, fixed at 200 here for reasonable computation time, is reached. 

The true conditional density is 

\[g^{0}\left(y_{1},y_{2}\right)=\pi_{1}^{0}f_{1}^{0}\left(y_{2}-F_{1}^{0}\left(y_{1}\right)\right)+\left(1-\pi_{1}^{0}\right)f_{2}^{0}\left(y_{2}-F_{2}^{0}\left(y_{1}\right)\right)\]

with $F_{i}^{0}\left(y_{1}\right)=a_{i}^{0}y_{1}+b_{i}^{0}$ and $f_{i}^{0}\sim\mathcal{N}\left(0,\left(\sigma_{i}^{0}\right)^{2}\right)$
for $i\in\left\{ 1,2\right\} $. 

For every example, we pick equal standard errors $\sigma_{1}^{0}=\sigma_{2}^{0}=0.5$, $b_{1}^{0}=0.5 $
and $b_{2}^{0}=-0.5 $ and let vary the rest of the coefficients: $\pi_{1}^{0}\in\left\{ 0.5,0.7,0.9\right\} $,
$a_{1}^{0},a_{2}^{0}\in\left\{ 0.1,0.5,0.9\right\} $. 

Let us focus in one particular example. Figure 1 illustrates one out of the twenty samples in the case 
$n=1000$, $\pi_{1}^{0}=0.7$, $a_{1}^{0}=0.1$ and $a_{2}^{0}=0.5$. The observed values of the series $Y_{t}$
are plotted on the first graph. On the second graph, we represent the convergence of the estimates for the mixture probabilities (solid and dashed lines) to the true values (dotted lines). Only the best result from the different 
initializations of the EM algorithm was drawn.

The summary of the results for all the examples is given by Table 1. In almost every case, the convergence is reached for the samples containing 2000 inputs. In practice, the results will be then more or less accurate, depending on the size of the sample, but also on the proximity of the components and on their frequency.

\begin{figure}[h!]
{\centering \resizebox*{10cm}{7cm}{\rotatebox{-90}{\includegraphics{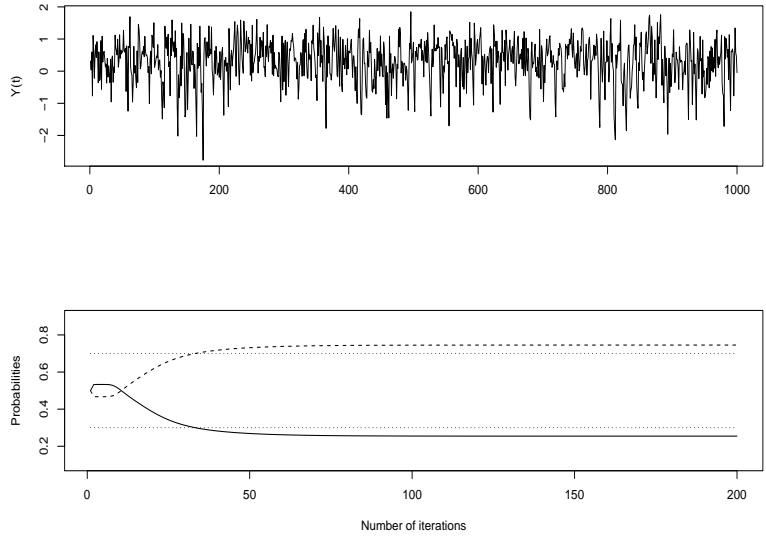}}} \par}
\caption{\small Series $Y_{t}$ and estimates for mixture probabilities $\pi_{i}^{0}$ \label{figure1}}
\end{figure} 

\begin{table}[h!]
\begin{center}\
\begin{tabular}{|cc|ccc|ccc|ccc|}
\hline
{\footnotesize $ $}&
{\footnotesize $\pi _{1}^{0}$}&
\multicolumn{3}{c|}{{\footnotesize 0.5}}&
\multicolumn{3}{c|}{{\footnotesize 0.7}}&
\multicolumn{3}{c|}{{\footnotesize 0.9}}\\
\cline{3-5} \cline{6-8} \cline{9-11}
{\footnotesize $ $}&
{\footnotesize $n$}&
{\footnotesize $\hat{p}=1$}&
{\footnotesize $\hat{p}=2$}&
{\footnotesize $\hat{p}=3$}&
{\footnotesize $\hat{p}=1$}&
{\footnotesize $\hat{p}=2$}&
{\footnotesize $\hat{p}=3$}&
{\footnotesize $\hat{p}=1$}&
{\footnotesize $\hat{p}=2$}&
{\footnotesize $\hat{p}=3$}\\
\hline
{\footnotesize $a_{1}^{0}=0.1$}&
{\footnotesize 200}&
{\footnotesize 20}&
{\footnotesize 0}&
{\footnotesize 0}&
{\footnotesize 20}&
{\footnotesize 0}&
{\footnotesize 0}&
{\footnotesize 20}&
{\footnotesize 0}&
{\footnotesize 0}\\
{\footnotesize $a_{2}^{0}=0.1$}&
{\footnotesize 500}&
{\footnotesize 18}&
{\footnotesize 2}&
{\footnotesize 0}&
{\footnotesize 18}&
{\footnotesize 2}&
{\footnotesize 0}&
{\footnotesize 20}&
{\footnotesize 0}&
{\footnotesize 0}\\
&
{\footnotesize 1000}&
{\footnotesize 14}&
{\footnotesize 6}&
{\footnotesize 0}&
{\footnotesize 9}&
{\footnotesize 11}&
{\footnotesize 0}&
{\footnotesize 11}&
{\footnotesize 9}&
{\footnotesize 0}\\
&
{\footnotesize 1500}&
{\footnotesize 6}&
{\footnotesize 14}&
{\footnotesize 0}&
{\footnotesize 4}&
{\footnotesize 16}&
{\footnotesize 0}&
{\footnotesize 5}&
{\footnotesize 15}&
{\footnotesize 0}\\
&
{\footnotesize 2000}&
{\footnotesize 5}&
{\footnotesize 15}&
{\footnotesize 0}&
{\footnotesize 0}&
{\footnotesize 20}&
{\footnotesize 0}&
{\footnotesize 1}&
{\footnotesize 19}&
{\footnotesize 0}\\
\hline 
{\footnotesize $a_{1}^{0}=0.1$}&
{\footnotesize 200}&
{\footnotesize 12}&
{\footnotesize 8}&
{\footnotesize 0}&
{\footnotesize 13}&
{\footnotesize 7}&
{\footnotesize 0}&
{\footnotesize 20}&
{\footnotesize 0}&
{\footnotesize 0}\\
{\footnotesize $a_{2}^{0}=0.5$}&
{\footnotesize 500}&
{\footnotesize 11}&
{\footnotesize 9}&
{\footnotesize 0}&
{\footnotesize 6}&
{\footnotesize 14}&
{\footnotesize 0}&
{\footnotesize 18}&
{\footnotesize 2}&
{\footnotesize 0}\\
&
{\footnotesize 1000}&
{\footnotesize 0}&
{\footnotesize 20}&
{\footnotesize 0}&
{\footnotesize 1}&
{\footnotesize 19}&
{\footnotesize 0}&
{\footnotesize 14}&
{\footnotesize 6}&
{\footnotesize 0}\\
&
{\footnotesize 1500}&
{\footnotesize 0}&
{\footnotesize 20}&
{\footnotesize 0}&
{\footnotesize 0}&
{\footnotesize 20}&
{\footnotesize 0}&
{\footnotesize 8}&
{\footnotesize 12}&
{\footnotesize 0}\\
&
{\footnotesize 2000}&
{\footnotesize 0}&
{\footnotesize 20}&
{\footnotesize 0}&
{\footnotesize 0}&
{\footnotesize 20}&
{\footnotesize 0}&
{\footnotesize 7}&
{\footnotesize 13}&
{\footnotesize 0}\\
\hline 
{\footnotesize $a_{1}^{0}=0.1$}&
{\footnotesize 200}&
{\footnotesize 0}&
{\footnotesize 20}&
{\footnotesize 0}&
{\footnotesize 4}&
{\footnotesize 16}&
{\footnotesize 0}&
{\footnotesize 17}&
{\footnotesize 3}&
{\footnotesize 0}\\
{\footnotesize $a_{2}^{0}=0.9$}&
{\footnotesize 500}&
{\footnotesize 0}&
{\footnotesize 20}&
{\footnotesize 0}&
{\footnotesize 0}&p
{\footnotesize 20}&
{\footnotesize 0}&
{\footnotesize 9}&
{\footnotesize 11}&
{\footnotesize 0}\\
&
{\footnotesize 1000}&
{\footnotesize 0}&
{\footnotesize 20}&
{\footnotesize 0}&
{\footnotesize 0}&
{\footnotesize 20}&
{\footnotesize 0}&
{\footnotesize 9}&
{\footnotesize 11}&
{\footnotesize 0}\\
&
{\footnotesize 1500}&p
{\footnotesize 0}&
{\footnotesize 20}&
{\footnotesize 0}&
{\footnotesize 0}&
{\footnotesize 20}&
{\footnotesize 0}&
{\footnotesize 4}&
{\footnotesize 16}&
{\footnotesize 0}\\
&
{\footnotesize 2000}&
{\footnotesize 0}&
{\footnotesize 20}&
{\footnotesize 0}&
{\footnotesize 0}&
{\footnotesize 20}&
{\footnotesize 0}&
{\footnotesize 0}&
{\footnotesize 20}&
{\footnotesize 0}\\
\hline
\end{tabular}
\end{center}
\caption{Number of components for $b_{1}^{0}=0.5$ and $b_{2}^{0}=-0.5$ }
\end{table}

\newpage 

\subsection{Laser time series}

A second example studies the complete laser series of ``Santa Fe time series prediction
and analysis competition''. The level of noise in this series is very low, the main source
of noise being the errors of measurement. We use the 12500 patterns for estimation. The Figure 2 shows the last 1000 patterns. The series is supposed to be stationary, and recall from Section 2 that a piecewise stationary time series is globally stationary if every component is stationary itself. The mixture of expert models is an example of piecewise stationary and globally stationary time series.

\begin{figure}[h]
{\centering \resizebox*{12cm}{3.5cm}{\rotatebox{-90}{\includegraphics{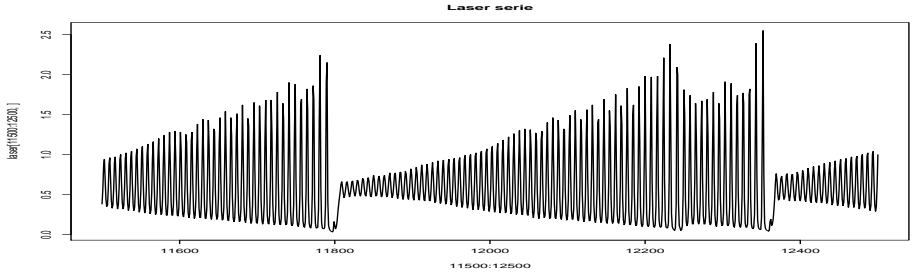}}} \par}
\caption{\small 1000 last Patterns of the laser series\label{figure2}}
\end{figure} 

We want to choose the number of components of the mixture by minimizing the BIC criteria.
Based on previous study (Rynkiewicz, 1999) we choose to use experts with 10 entries, 5 hidden units, one linear output, and hyperbolic tangent activation functions. We want to know the number of experts to use with such series. As this is a real application, it is impossible to check the main assumption of our theory : the true model belongs to the set of possible models. However, we want to know if the developed theory can give an insight for choosing the number of experts.

The parameters are estimated using the standard EM algorithm.  In order to avoid bad local maxima, estimation is performed with 100 different initializations of model parameters. For each estimation we proceed  with 200 iterations of the EM algorithm and for each M-step we optimize the parameter of the MLPs until their error of prediction doesn't improve anymore.  

The estimated loglikelihood and the estimated probability of the mixture are the following : 

\begin{center}
\begin{tabular}{|c|c|c|}
\hline
number of experts&BIC&probabilities of the mixture\\
\hline
1&-32.16894&1\\
2&-25.92&(0.8,0.2)\\
3&-38.42&(0.7,0.21,0.09)\\
\hline
\end{tabular}
\end{center}

The results are clear for our model, the best model is the model with two experts. It is difficult to give an interpretation of the regimes because mixing probabilities remain constant over time. However, if we look at the prediction made by each expert, we  can see that one expert seems to be specialized in the general regime of the series and the second one with the collapse regime.

\bigskip{}

The proposed method gives an insight on the way to choose the number of experts in a mixture model for laser time series. However, since the probabilities of the mixture are constant, it would be better to choose probabilities depending on the previous value of the time series as in the gating expert of Weigend et al. (Weigend et. al, 1995) or of the time as in hybrid hidden Markov/MLP Models (Rynkiewicz, 2006). The prediction error of this simple mixture model is not competitive with such more complex models, however we need to improve the theory to deal with such complex modeling.  

\section{Conclusion and future work}

We have proven the consistency of the BIC criterion for estimating the number of components in a mixture of multilayer perceptrons. In our opinion, two important directions are to be studied in the future. The case of mixtures should be extended to the general case of gated experts which allow the probability distribution of the multilayer perceptrons to depend on the input and thus, to learn how to split the input space. The second possible extension should remove the hypothesis of a fixed number of units on the hidden layer. The problem of estimating the number of hidden units in one multilayer perceptron was solved in Rynkiewicz (2006), but it would be interesting to mix the two results and prove the consistency of a penalized criterion when there is a double non-identifiability problem : number of experts and number of hidden units.

\footnotesize{

}

\end{document}